\documentclass[11pt]{article}

\usepackage{graphicx}
\usepackage{amsmath}
\usepackage{amsthm}
\usepackage{amssymb}
\usepackage[mathscr]{eucal}

\title{The intrinsic asymmetry and inhomogeneity\\ of Teichm\"{u}ller space}
\author{Benson Farb and Shmuel Weinberger \thanks{Both authors are
supported in part by the NSF.}}

\theoremstyle{plain}
\newtheorem{theorem}{Theorem}[section]

\def\endproof{$\diamond$ \bigskip}

\def\title{\em}

\newcommand\R{\mbox{\bf R}}

\newcommand\Z{\mbox{\bf Z}}

\newcommand\Q{\mbox{\bf Q}}

\DeclareMathOperator\vcd{vcd}
\DeclareMathOperator\Fix{Fix}

\DeclareMathOperator\Out{Out}

\DeclareMathOperator\teich{Teich}
\DeclareMathOperator\Teich{\teich}

\DeclareMathOperator\Mod{Mod}

\DeclareMathOperator\M{{\cal M}}

\DeclareMathOperator\genus{genus}

\DeclareMathOperator\Isom{Isom}

\DeclareMathOperator\QD{QD}

\DeclareMathOperator\sol{sol}
\DeclareMathOperator\sem{ss}

\begin{document}
\maketitle


\section{Introduction}

Let $S=S_{g,n}$ be a connected, orientable 
surface of genus $g\geq 0$ with $n\geq 0$ punctures.  Let $\Teich(S)$ denote the
corresponding Teichm\"{u}ller space, and let $\Mod(S)$ denote the
mapping class group of $S$.  Understanding the analogy of
$\Teich(S)$ with symmetric spaces is a well-known theme.  Recall
that a complete Riemannian manifold $X$ is {\em symmetric} if it
is symmetric at each point $x\in X$: the map $\gamma(t)\mapsto
\gamma(-t)$ which flips
geodesics about $x$ is an isometry.  Symmetric spaces $X$ are 
{\em homogeneous}: the isometry group $\Isom(X)$ acts
transitively on $X$.  In his famous paper \cite{Ro}, 
Royden studied the possible symmetry and homogeneity of
$\Teich(S)$, endowed with the Teichm\"{u}ller metric.

\begin{theorem}[Royden \cite{Ro}]
\label{theorem:Royden}
Suppose $S$ is closed of genus at least $2$, and let $\Teich(S)$
be Teichm\"{u}ller space endowed with the Teichm\"{u}ller
metric $d_{Teich}$. Then
\begin{enumerate}
\item[a.] $(\Teich(S),d_{\Teich})$ is not symmetric at any point.
\item[b.] $\Isom(\Teich(S),d_{\Teich})$ contains $\Mod(S)$ (modulo its center if $S$ is closed of 
genus $2$) as a subgroup of
index $2$.
\end{enumerate}
\end{theorem}

Note that the conclusion of Theorem \ref{theorem:Royden} is false
when $\genus(S)=1$, as $d_{\Teich}$ in this case is the
hyperbolic metric on the upper half-plane.  Earle-Kra \cite{EK} extended part (b) of 
Royden's Theorem 
to arbitrary surfaces of finite type.  Royden deduced
Theorem \ref{theorem:Royden} from a detailed analysis of the fine
structure of the space $\QD^1(M)$ of unit norm holomorphic
quadratic differentials on a Riemann surface $M$.  In particular,
he found an embedding of $M$ in $\QD^1(M)$ and characterized it
by the degree of Holder regularity of the norm on $\QD^1(M)$ at
the points of the embedding.

The Teichm\"{u}ller metric is a complete Finsler metric, under
which moduli space ${\cal M}(S):=\Teich(S)/\Mod(S)$ has finite
volume (see the proof of Theorem 8.1 of \cite{Mc}) \footnote{A
Finsler metric determines a unique volume form by declaring the
unit Finsler ball at each point to have volume $1$.}.  There are
many other complete, finite volume, $\Mod(S)$-invariant Finsler
(indeed Riemannian) metrics on $\Teich(S)$, each with special
properties.  Examples include the Kahler-Einstein metric,
McMullen's metric, and the (perturbed) Ricci metric.  Following
Royden's approach to understanding the symmetry and homogeneity of
these metrics appears to involve difficult analysis.

The goal of this paper is to explain a completely different, nonanalytic mechanism behind 
Royden's Theorem.  It will allow us to extend much of his result from the Teichm\"{u}ller metric to 
any metric, including each of those mentioned above.  The theorem is the following.

\begin{theorem}
\label{theorem:main}
Let $S=S_{g,n}$ be a surface, and let 
$d$ be any complete, finite covolume, $\Mod(S)$-invariant Finsler (e.g. Riemannian) 
metric on $\Teich(S)$.  Then
\begin{enumerate}
\item[a.] If $g\geq 3$ then $(\Teich(S),d)$ is not symmetric at any point.
\item[b.] If $3g-3+n\geq 2$ then $\Isom(\Teich(S),d)$ contains $\Mod(S)$  (modulo its center if $(g,n)=(1,2)$ or $(2,0)$) as a subgroup of finite index.
\end{enumerate}
\end{theorem}

In this way Teichm\"{u}ller space exhibits a kind of intrinsic asymmetry and inhomogeneity.  The number $3g-3+n$ plays the role of  ``$\Q$-rank'' in this context.  Theorem \ref{theorem:main}(b) is sharp: the conclusion is false whenever $3g-3+n<2$, and indeed the corresponding Teichm\"{u}ller spaces admit 
hyperbolic metrics.   

\medskip
\noindent
{\bf Remarks. }
\begin{enumerate}

\item The proof of Theorem \ref{theorem:main} gives immediately 
that the hypothesis can be weakened to allow $d$ 
to be any $\Gamma$-invariant metric, for any finite index subgroup $\Gamma<\Mod(S)$.  As such subgroups are typically torsion free, this shows that the phenomenon of inhomogeneity and asymmetry is not the result of the constraints imposed by having finite order symmetries of the metric.

\item Ivanov proved (see, e.g., \cite{I2}) that the moduli space of
Riemann surfaces $\Teich(S)/\Mod(S)$ never admits a locally
symmetric metric when $S$ is closed and $\genus(S)\geq 2$.  If it
did admit such a metric, then $\Teich(S)$ would admit a complete, finite covolume, 
$\Mod(S)$-invariant metric which is symmetric at {\em every} point.
Thus Theorem \ref{theorem:main}(a) gives a new proof, and 
generalization, of Ivanov's theorem.
\end{enumerate}

\smallskip
The proof of Theorem \ref{theorem:main}(a) relies on Theorem \ref{theorem:main}(b).  One key  ingredient in our proofs is 
Smith theory.  This is not the first time Smith theory has been used to analyze 
actions on Teichm\"{u}ller space: Fenchel used it in the 1940's
to analyze certain periodic mapping classes.

We conjecture that the index in Theorem \ref{theorem:main}(b) can be taken to depend only on $S$.  More strongly, one might hope that it can always be taken to be $2$.  As evidence towards this 
strongest possible conjecture, we can prove it in the ``$\Q$-rank $2$'' case.

\begin{theorem}
\label{theorem:exact}
Let $S$ be the twice-punctured torus or the $5$-punctured sphere.   Then 
the index in Theorem \ref{theorem:main}(b) can be taken to be $2$.
\end{theorem}

\medskip
\noindent
{\bf Acknowledgements. }We are grateful to Tom Church, Alex Eskin and Dan Margalit for their 
numerous useful and insightful comments.  We are especially grateful to Yair Minsky, whose persistence in pinning us down, and whose many comments and corrections, 
greatly improved this paper.

\section{Proof of Theorem \ref{theorem:main}(b)}
\label{section:royden}

By the Myers-Steenrod Theorem (or, in the general Finsler case,
Deng-Hou \cite{DH}), the group $I:=\Isom(\Teich(S),d)$ is a Lie group, possibly
with infinitely many components, acting properly discontinuously on $\Teich(S)$.  We remark 
that Theorem 2.2 
of \cite{DH} states that any isometry of a Finsler metric is necessarily a diffeomorphism; we will 
use this smoothness later.   

Let $I_0$ denote the connected
component of the identity of $I$; note that $I_0$ is normal in $I$.   Let $\mu$ denote the 
measure on $\Teich(S)$ (and the induced measure on $\Teich(S)/\Mod(S)$) induced by the 
Finsler metric $d$.  We are assuming that $\mu(\Teich(S)/\Mod(S))<\infty$.

If $I$ is discrete, then we claim that $[I:\Mod(S)]<\infty$.  To see this, let $K_I$ be a measure-theoretic fundamental domain for the $I$-action on $\Teich(S)$.  By this we mean that:
\begin{enumerate}
\item The complement of $\bigcup_{g\in I}g\cdot K_I$ in $\Teich(S)$ has $\mu$-measure $0$, and
\item For all nontrivial $g\in I$, we have $\mu(g\cdot K_I\cap K_I)=0$.
\end{enumerate}

Note that measurable fundamental domains always exist for properly discontinuous actions.  
Let $K_{\Mod(S)}$ be a measure-theoretic fundamental domain for $\Mod(S)$.  Let $\{g_i\}$ be a collection of coset representatives for $\Mod(S)$ in $I$.  Then we can choose $K_{\Mod(S)}$ so that 
$$K_{\Mod(S)}=\bigcup_{i}g_i\cdot K_I.$$
Notice that this union is a disjoint union, up to sets of measure zero, by Condition (2) in the definition of fundamental domain.   Since by assumption $\mu(K_I)$ and $\mu(K_{\Mod(S)})$ are both finite, it follows that this must be a finite union, i.e. $[I:\Mod(S)]<\infty$.  

So suppose that $I$ is not discrete.  Myers-Steenrod then gives that the dimension of $I$ is
positive, and so $I_0$ is a connected, 
positive-dimensional Lie group.   Let $\Gamma:=\Mod(S)$, and let $\Gamma_0:=I_0\cap \Mod(S)$.   We have the following exact sequences:

\begin{equation}
\label{eq:exact1}
1\to I_0\to I\to I/I_0\to 1
\end{equation}

and 
\begin{equation}
\label{eq:exact2}
1\to \Gamma_0\to \Gamma \to \Gamma/\Gamma_0\to 1
\end{equation}

\medskip \noindent {\bf Step 1 ($\Gamma_0$ is a lattice in
$I_0$): }   We begin by replacing $\Gamma$ (hence $\Gamma_0$) by a torsion-free subgroup of finite index, but keep the notation $\Gamma$ (resp. $\Gamma_0$) for this new group.   While not formally necessary, this assumption will make ``orbifold technicalities'' disappear; the point is that quotients of $\Teich(S)$ and associated bundles by $\Gamma$ will be manifolds.

We begin by replacing the proper action of $I$ on $\Teich(S)$ by a free action on an associated space ${\cal F}(\Teich(S))$, which we now describe.

For any smooth, connected, $n$-dimensional Finsler manifold $M$ 
we have the natural unit sphere-bundle over $M$, which is a sub-bundle of the tangent bundle $TM$, with fiber over $m\in M$ the unit Finsler sphere $S_m$ of $TM_m$.  We also have the associated bundle $E\to M$ whose fiber is the $n$-fold product $S_m^n$.  Let ${\cal F}(M)$ denote the sub-bundle of this bundle with fiber the set of $n$-tuples of distinct points of $S_m$ that span $TM_m$.   The group $\Isom(M)$ clearly acts on ${\cal F}(M)$ by homeomorphisms.  

The exponential map on a Finsler manifold is a local diffeomorphism (see, e.g.\, \cite{DH}, Lemma 1.1).  Since we also have that Finsler isometries take geodesic rays to geodesic rays, we see that the set of points of $M$ for which a Finsler isometry is the identity and has derivative the identity is both open and closed.  Thus the action of $\Isom(M)$ on ${\cal F}(M)$ is free.

Now, we will want to apply the Slice Theorem to this action (see below), and to do so we need it to preserve some smooth structure on ${\cal F}(M)$.  To this end, first note that $M$ is assumed to be a smooth manifold, and since $\Isom(M)$ acts properly on $M$, we have that $\Isom(M)$ actually preserves some Riemannian metric on $M$.  We then note that there is an $\Isom(M)$-equivariant homeomorphism from ${\cal F}(M)$ to the bundle of unit $n$-frames over $M$, and so the 
action of $\Isom(M)$ on ${\cal F}(M)$ preserves the pullback of the standard smooth structure on the latter bundle.  We remark that the reason we do not simply replace the given Finsler metric with this invariant Riemannian metric is that the finiteness of volume of $M$ in the Finsler metric need not imply such finiteness for the invariant Riemannian metric.

We now wish to construct an $\Isom(M)$-invariant measure on ${\cal F}(M)$.  To this end, note that 
the bundle $E\to M$ discussed above is locally a product $U\times S^n$, where $U$ is a neighborhood in $M$.  The Finsler metric on $M$ determines a volume form on $M$, which induces a measure $\nu$ on $U$.  On $S$, we have an induced measure $\mu$ which is given infinitessimally by the rule that, for a subset $A\subseteq S_m$, the measure is given by the measure (induced by the Finsler norm of $A$ on $TM_m$) of the Euclidean cone of $A$, normalized so that the measure of $S_m$ equals $1$.  The local product measure $\nu\times \mu$ then gives an $\Isom(M)$-invariant measure on $E$, which in turn induces an $\Isom(M)$-invariant 
measure on ${\cal F}(M)$.  By construction, the pushforward of this measure under the natural projection ${\cal F}(M)\to M$ is the measure on $M$ induced by the given Finsler metric; in particular, if $M$ is assumed to have finite measure then ${\cal F}(M)$ has finite measure.

Now let $M=\Teich(S)$.  Pick $x\in{\cal F}(M)$, and consider the $I$-orbit $O_x$.    The Slice Theorem for proper group actions (see, e.g., \cite{DK}, 2.4.1) 
in this context asserts that there is an $I$-invariant tubular
neighborhood $V$ of $O_x$ in ${\cal F}(\Teich(S))$ that is a homogeneous 
vector bundle $\pi:V\to O_x$.  The measure on ${\cal F}(M)$ constructed above restricts to a measure on $V$, and $\pi$ pushes this forward to a left-invariant measure on $O_x$, which we can identify with a left-invariant measure on $I$.   Note that all left-invariant measures on $I$ 
are proportional, by uniqueness of Haar measure.  The key property we will use is that if a subset $A\subseteq I$ has infinite measure then $\pi^{-1}(A)$ has infinite measure.

Choose any fiber $D$ of the bundle $V\to O_x$, so that $V$ is the $I$-orbit of $D$.   Since $I_0$ 
is the connected component of the identity of $I$, and since $I_0$ is a closed subgroup of $I$, we have that $V$ can be written as a disjoint union of $I_0$-orbits of $D$, one for each element of $\pi_0(I)$.  Note that $\Gamma/\Gamma_0 \subseteq \pi_0(I)$.  Thus $V/\Gamma$ 
is given by the image of the $I_0$-orbit $W$ of $D$ under the projection 
$${\cal F}(\Teich(S))\to {\cal F}(\Teich(S))/\Gamma=(\Teich(S)/\Gamma)={\cal F}(\M(S)).$$

Since the $I$-action on ${\cal F}(\Teich(S))$ is free, this projection is a measure-preserving homeomorphism when restricted to $W$.  Now if $I_0/\Gamma_0$ had infinite measure, then so would $W$ (by the discussion above), and thus so would ${\cal F}(\M(S))$.  On the other hand the map ${\cal F}(\M(S))\to \M(S)$ is measure-preserving by the construction of the measure on ${\cal F}(\M(S))$ (see above).  This gives that $\M(S)$ has infinite measure, contradicting the given.  We conclude that $I_0/\Gamma_0$ has finite measure, as desired.

\medskip 
\noindent 
{\bf Step 2 ($I_0$ is semisimple with finite
center): } 
For any connected Lie group $G$ there is an exact sequence 
\begin{equation}
\label{equation:gen1}
1\to G^{\sol}\to G\to G^{\sem} \to 1
\end{equation}
where $G^{\sol}$ denotes the {\em solvable radical} of $G$ (i.e. the maximal
connected, normal, solvable Lie subgroup of $G$), and where $G^{\sem}$ is
the connected semisimple Lie group $G/G^{\sol}$.  

We can apply this setup with $G=I_0$.   We now follow exactly the argument of the proof of Proposition 
3.3 of \cite{FW}, with the only change being that here $\Gamma$ is not assumed torsion free.  First, a theorem of Raghunathan gives a unique maximal normal solvable subgroup $\Gamma_0^{\sol}$ of $\Gamma_0$.  Since $\Gamma_0^{\sol}$ is unique it is characteristic in $\Gamma_0$, hence normal in $\Gamma$.  The main theorem of 
\cite{BLM} states that any solvable subgroup of $\Gamma=\Mod(S)$ has an abelian subgroup of finite index.  Hence $\Gamma_0^{\sol}$ has a torsion-free, normal (in $\Gamma$) finite-index abelian subgroup, which we denote by $N$.  We will prove below that $\Gamma=\Mod(S)$ has no infinite normal abelian subgroups, from which we will conclude that $N$ is trivial, so that $\Gamma_0^{\sol}$ is finite.

Now the next part of the argument of Proposition 3.3 of \cite{FW} quotes a theorem of Prasad.  In applying Prasad's theorem, we used the fact that the sum of the ranks of the abelian quotients of the derived series of $\Gamma_0^{\sol}$ equals $0$.  While in \cite{FW} this followed from the fact that $\Gamma^{\sol}_0=0$, we only need that $\Gamma_0^{\sol}$ is finite, which we have just proven.   We may thus quote Prasad's theorem and conclude, as in \cite{FW}, that $I_0^{\sol}$ is compact, and that the center $Z(I_0)$ is finite.  

Given this, we finish the proof of Step 2 as follows.    Since any compact Lie group is the product of a simple Lie group and a torus, $I_0^{\sol}$ must be a torus.  Further, as is precisely argued in 
Proposition 3.3 of \cite{FW},  the conjugation action of the connected group 
$I_0^{\sem}$ on the torus $I_0^{\sol}$, which has discrete automorphism group, must be trivial, so that $I_0^{\sol}$ must be a direct factor of $I_0$.  We thus need to rule this out when $I_0^{\sol}$ is positive-dimensional.  

So if we can prove that $\Mod(S)$ has no infinite, normal
abelian subgroup $A$, and if we can rule out that $I_0^{\sol}$ is a torus direct factor of $I_0$, then we have completed Step 2.  We begin with the first claim.  

By the classification of abelian subgroups
of $\Mod(S)$ (see \cite{BLM} or \cite{I1}), any abelian subgroup
$A$, after perhaps being replaced by a finite index
characteristic subgroup if necessary, 
either is cyclic with a 
pseudo-Anosov generator or there is a unique maximal 
finite collection $C$ of simple
closed curves, called the {\em canonical reduction system} 
of $A$, left invariant (setwise) by each $a\in A$.  In the
first case, the normalizer of $A$ is virtually cyclic, a
contradiction, so suppose we are in the latter case.  Then for
any $f\in\Mod(S)$ the canonical reduction system for $fAf^{-1}$
is $f(C)$.  The result now follows by picking an $f$ such that
$f(C)\neq C$.  Thus $A$ must be trivial.

We now rule out that $I_0$ contains a
positive-dimensional torus as a direct factor.  Suppose it does.   Then $Z(I_0)$ is positive-dimensional.  Since $Z(I_0)$ is a positive-dimensional abelian Lie group, and we can write $Z(I_0)=A\times (S^1)^d\times \R^k$ for some $d>0,k\geq 0$.  Now let $T$ denote the maximal compact subgroup of the connected component of the identity of $Z(I_0)$.  From the above description of $Z(I_0)$, it is clear that $T$ is characteristic in $Z(I_0)$, hence in $I_0$.  
Since $I_0$ is normal in $I$, we have
that the conjugation action of $\Gamma=\Mod(S)$ leaves $T$
invariant.  Thus the action of $\Gamma$ on $\Teich(S)$ leaves
$\Fix(T)$ invariant.  Since $T$ is acting smoothly,
$\Fix(T)$ is a manifold.  Since $\Teich(S)$ is contractible, $\Fix(T)$ is
acyclic (see, e.g., \cite{Br}, Theorem 10.3).  As $T$ is positive
dimensional and connected, and since the $T$-action on
$\Teich(S)$ is faithful, $\dim(\Fix(T))<\dim(\Teich(S))-1$.  Note
that the action of $\Mod(S)$ on $\Fix(T)$ is properly
discontinuous, being the restriction of the properly
discontinuous action of $\Mod(S)$ on $\Teich(S)$.

\bigskip \noindent {\bf Case A ($\dim(\Fix(T))>3$): } We claim
that there is a contractible manifold $Z$ of dimension
$\dim(\Fix(T))+1<\dim(\Teich(S))$ on which $\Mod(S)$ acts
properly discontinuously.  Given this, we recall that Despotovic 
\cite{D} proved that $\Mod(S)$ admits no properly
discontinuous action on any contractible manifold of dimension
$<\dim(\Teich(S))$, giving us a contradiction.  Thus $T$ would be
trivial, and so $I_0$ is semisimple with finite center.

We now prove the claim.    Note that by Smith theory (see \cite{Br}) $\Fix(T)$ is acyclic.  We want
to replace this $\Mod(S)$-manifold by another one that is contractible, i.e.\ we want 
to kill the fundamental group.  If $\Fix(T)/\Mod(S)$ is compact, then
the kernel that we are trying to kill is finitely normally
generated, and the construction is standard:  one kills the elements
of this kernel by surgering the circles, giving rise to new homology =
homotopy in dimension $2$, which can then be killed by surgering the
$2$-spheres as well (see Section 3 of [H]).  Taking the product of $\Fix(T)$ with $\R$ if 
$\dim(\Fix(T))=4$,  the $2$-spheres needed at this stage can be embedded by
general position since $\dim(\Fix(T)\times \R) > 4$.

If $\Fix(T)/\Mod(S)$ is not compact, then one wants to do the same
argument, but one has to be careful to make sure that the circles and
$2$-spheres that one wants to surger do not accumulate.   
However, by replacing $\Fix(T)$ by $\Fix(T)\times \R$, this problem disappears: 
one simply does the $i^{\rm th}$ surgery at the ``height'' $\Fix(T)\times \{i\}$, 
producing at height i a 2-dimensional homology class.  This class
can be represented by a sphere which lies in a compact region that is
above level $i-1/2$ (by the Hurewicz theorem) that is embedded (by
general position) with trivial normal bundle (by appropriate
framing).  The infinite collection of constructed $2$-spheres do not
accumulate, since at most $i$ of them pass through level $i$.  Surgering
this free basis for the homology of the previous stage gives us back
an acyclic manifold, which is now, in addition, simply connected, and,
thus, contractible.

\bigskip
\noindent
{\bf Case B ($\dim(\Fix(T))\leq 3$): }In this case the proper action of $\Mod(S)$ on the acyclic space 
$\Fix(T)$ implies that the virtual cohomological dimension
$\vcd(\Mod(S))$ satisfies 
$$\vcd(\Mod(S))\leq \dim(\Fix(T))\leq 3.$$ 
But  then by the formulas for $\vcd(\Mod(S))$, given for example in Theorem 6.4 of \cite{I2}, and 
since $3g-3+n\geq 2$ by hypothesis, this leaves the cases of possible $(g,n)$ to be 
one of $\{(2,0), (0,5), (0,6), (1,2), (1,3)\}$.   As the action of $\Mod(S)$ on $\Fix(T)$ is properly discontinuous, this rules out 
$\dim(\Fix(T))=0$.   Now $\Fix(T)$ has even codimension in 
$\Teich(S)$ (see \cite{Br}, Theorem 10.3), and so is
even-dimensional.  As we are assuming in this case that $\dim(\Fix(T))\leq 3$, it follows 
that $\dim(\Fix(T))=2$.  But this would imply that $\Mod(S)$ has a 
a (closed or open) surface group as a subgroup of finite index.  It cannot, however, because for instance in these cases $\Mod(S)$ contains both $\Z^2$, eliminating all 
surfaces but the torus,  and also 
a rank $2$ free group, eliminating the torus.

\medskip 
\noindent 
{\bf Step 3 ($I_0$ has no compact factors): }  Let $K$ be the maximal compact factor of $I_0$.  Since 
$I_0$ is semisimple with finite center, $K$ is characteristic in $I_0$.  Since $I_0$ is normal in 
$I$, it follows that $K$ is invariant under conjugation by any element of $I$.  Further, note that 
$K$ is semisimple since $I_0$ is semisimple.

Since $K$ is compact, we have
that $K\cap\Gamma=K\cap\Gamma_0$ is a finite normal subgroup of $\Gamma$.   Replacing $\Gamma$ by a finite index subgroup, which we will also denote by $\Gamma$, we can assume that $K\cap \Gamma$ is trivial.  As $K$ is invariant by conjugation by elements of $I$, we have an exact sequence 
\begin{equation}
\label{eq:1}
1\to K\to \langle K,\Gamma\rangle\to \Gamma\to 1
\end{equation}
where the middle term denotes the subgroup of $I$ generated by $K$ and $\Gamma$.  
As explained in IV.6 of \cite{Bro}, any exact sequence 
$$1\to A\to B\to C\to 1$$ is determined by two pieces of data: a representation 
$\rho:C\to \Out(A)$ and a cocycle $\eta\in H^2(C,Z(A))$, where $Z(A)$ denotes the center of $A$, and is a $C$-module via the action of $\rho$.  In the case (\ref{eq:1}), we have that both $Z(K)$ and $\Out(K)$ are finite since $K$ is semisimple.   
We may thus pass to a finite index
subgroup $\Lambda<\Gamma$ so that $\rho$ has trivial image.  Now let
$$1\to Z(K)\to \widehat{\Lambda}\to\Lambda\to 1$$
be the group extension corresponding to 
the pullback of the cocycle $\eta\in H^2(\Lambda,Z(K))$ 
corresponding to the extension (\ref{eq:1}) restricted to $\Lambda$.   Note that this 
exact sequence defines a trivial cocycle.  Note also that $\widehat{\Lambda}\subseteq \langle K,\Lambda\rangle$.  Now $Z(K)$ is central both 
in $\widehat{\Lambda}$ and in $K$.  We thus have a split exact sequence
\begin{equation}
\label{eq:2}
1\to K/Z(K)\to \langle K, \widehat{\Lambda}\rangle \to\widehat{\Lambda}\to 1
\end{equation}
Since $\rho$ has trivial image, we can change the section of (\ref{eq:2}) to get a copy of 
$K/Z(K)\times \widehat{\Lambda}$ in $\langle K, \widehat{\Lambda}\rangle $.   As noted above, this group is a subgroup of $\langle K,\Gamma\rangle$, and so acts on $\Teich(S)$.

If $K$ is positive-dimensional then so is 
$K/Z(K)$, and so $K/Z(K)$ contains a closed subgroup isomorphic to a circle $T$.  
As the (possibly noneffective) action of 
$\widehat{\Lambda}$ on $\Teich(S)$
commutes with the action of $T$, we have that $\widehat{\Lambda}$
leaves $\Fix(T)$ invariant.  But this gives a contradiction,
exactly as in Step 2 above, once we observe that \cite{D} applies to $\widehat{\Lambda}$, 
and so we obtain that $K$ is
trivial.  

To see that \cite{D} applies to $\widehat{\Lambda}$, there are two minor issues: her result is stated for $\Gamma=\Mod(S)$, while we need the theorem for finite extensions and 
finite index subgroups of $\Gamma$.  For the first issue we simply note that,
just as mentioned in the first sentence of the proof of
Theorem 26 in \cite{BKK}, the groups for which the theorem in
\cite{D} holds are closed under finite extensions.
The proof  that \cite{D} holds not just for $\Mod(S)$,
but for any finite index subgroup $\Gamma'$ of $\Mod(S)$, is 
verbatim the same as in \cite{D}, replacing the ``Mess
subgroups'' $B_g$ constructed there with $B_g\cap\Gamma'$, which has finite index in $B_g$.   The 
two key properties of $B_g$ are:
\begin{itemize}
\item $B_g$ is the fundamental group of a closed, triangulable topological manifold of dimension $4g-5$, and
\item The natural ``point pushing'' subgroup of $B_g$ is a closed surface group.
\end{itemize}

Each of these properties is clearly preserved by taking finite index subgroups, and so the proof of the main theorem of \cite{D} goes through when $\Mod(S)$ replaced by any finite index subgroup $\Gamma'$.

\medskip 
\noindent 

{\bf Step 4 ($I_0$ is trivial): } By the
previous steps, we know that $\Gamma_0$ is a lattice in the
semisimple Lie group $I_0$, and that $I_0$ has no compact
factors.  The proof of Proposition 3.1 of \cite{FW}
now gives that there is a finite index subgroup $\Gamma'$ of
$\Gamma$ so that 
\begin{equation}
\label{eq:fwprop}
\Gamma' \approx \Gamma_0 \times \Gamma'/
\Gamma_0.
\end{equation}
To give an idea of the proof of (\ref{eq:fwprop}) from \cite{FW}, we begin by considering the exact sequence 
\begin{equation}
\label{eq:exact5}
1\to \Gamma_0\to \Gamma\to \Gamma/\Gamma_0\to 1.
\end{equation}
As mentioned above, the extension, (\ref{eq:exact5}) is
determined by a representation $\rho:\Gamma/\Gamma_0\to \Out(\Gamma_0)$, and by a 
cohomology class in $H^2(\Gamma/\Gamma_0,Z(\Gamma_0)_\rho)$.  

One shows that the image of $\rho$ actually lies in $\Out(I_0)$, which is finite since 
$I_0$ is semisimple with finite center.   After replacing $\Gamma$ by a finite index subgroup $\Gamma'$, one then gets that the resulting representation $\rho$ is trivial.  One can also choose $\Gamma'$ so that $\Gamma_0\cap\Gamma'$ is torsion-free.  Since this group is a lattice in the semisimple Lie group with finite center $I_0$, it has finite center; since $\Gamma_0\cap\Gamma'$ is torsion-free, its center is trivial, 
so that the pertinent $H^2$ vanishes.

We now claim that for any finite index subgroup $\Gamma'<\Mod(S)$, if $\Gamma'=A\times B$ then 
either $A$ or $B$ is finite.  To see this, note that any such $\Gamma'$ contains a pseudo-Anosov homeomorphism $f$ (for example take a sufficiently high power of any pseudo-Anosov in $\Mod(S)$).  
The centralizer  in $\Gamma'$ (indeed in $\Mod(S)$) 
of any power of $f$ has $\Z$ as a finite index subgroup (see \cite{I1}, Lemma 8.13).  But in a product 
of two infinite groups, it is easy to see that any element has some power whose centralizer 
does not contain $\Z$ as a finite index subgroup.

Thus either $\Gamma_0$ is finite or $\Gamma'/\Gamma_0$ is finite.  The latter possibility 
implies that $\Gamma'$ , hence $\Mod(S)$, has a finite index subgroup which is isomorphic to a lattice 
in the semisimple Lie group $I_0$.  If $I_0$ is nontrivial, then it must contain a noncompact 
factor (by Step 3).  This would then contradict the theorem of Ivanov (see, e.g., \S 9.2 of \cite{I2}) 
that no finite index subgroup of $\Mod(S)$ is isomorphic to a lattice in 
a noncompact semisimple Lie group.  Thus it must be that either $I_0$ is trivial, or $\Gamma_0$ is 
finite.  If the latter possibility were to occur, then $I_0$ would be compact since $\Gamma_0$ is a lattice in $I_0$ by Step 1.  But this would contradict Step 3.
\endproof

\section{Proof of Theorem \ref{theorem:main}(a)}

Let $\tau$ be a {\em symmetry} of $(\Teich(S),d)$, i.e.\ an isometric involution with an isolated fixed point.  Let $L=\langle\Mod(S),\tau\rangle$ be the group generated by $\Mod(S)$ and
by $\tau$.

By Theorem \ref{theorem:main}(b), which we have already proven, $[L:\Mod(S)]<\infty$.  
Thus the action of $\tau$ on $L$ by conjugation induces a {\em
commensuration} of $\Mod(S)$, i.e.\ an isomorphism between two
finite index subgroups.  Since $\Mod(S)$ is residually finite, we
can pass to further finite index subgroups so that neither
contains the hyperelliptic involution.  By a theorem of Ivanov 
(see Theorem 8.5A of \cite{I2}), since $\genus(S)\geq 2$ any such commensuration agrees on
some finite index characteristic subgroup $H$ of $\Mod(S)$ with conjugation by some
element $\phi$ of the {\em extended mapping class group} $\Mod^\pm(S)$,
the index $2$ supergroup of $\Mod(S)$ which includes an
orientation-reversing homotopy class of homeomorphism.

We now claim that there exists an infinite order element $\psi_2\in \Mod(S)$ that
commutes with $\tau$.  Note that since the conjugation action of $\tau$ on $H$ agrees with 
the conjugation action of $\phi$,  it is enough to produce an infinite order 
element $\psi_2\in H$ so that $\psi_2$ commutes with $\phi$.

Given this claim, we complete the proof of the theorem as follows.  We are given that 
$\tau$ has an isolated fixed point $x\in\Teich(S)$.  By Smith theory, $\Fix(\tau)$ is $\Z/2\Z$ acyclic; in particular $\Fix(\tau)$ is connected.  Since $\psi_2$ is an infinite order mapping class, we 
have that $\psi_2(x)\neq x$, by proper discontinuity of the action of $\Mod(S)$ on $\Teich(S)$.  But 
$$\tau(\psi_2(x))=\psi_2(\tau(x))=\psi_2(x)$$
so that $\tau$ also fixes $\psi_2(x)\neq x$.  As $\Fix(\tau)$ is 
connected and fixes at least two distinct points, it must have positive 
dimension.  This contradicts the fact that $x$ is an isolated fixed point of $\tau$.  Thus such a 
$\tau$ cannot exist, and we are done.

We now prove the claim.  First note that 
since $\tau^2={\rm Id}$, conjugation by $\phi^2$
is the identity on some finite index subgroup $H$ of
$\Mod(S)$. Now there exists $N>0$ so that for a Dehn twist $T_\alpha$ about any simple 
closed curve $\alpha$, we have  $T^{N}_\alpha\in H$.  For any twist $T_\alpha$
and any element $f\in\Mod(S)$, we have the well-known formula
$$fT^N_\alpha f^{-1}=T^N_{f(\alpha).}$$
Since $\phi^2\in\Mod(S)$, we can apply this formula 
with $f=\phi^2$, giving that $T^N_\alpha=T^N_{\phi^2(\alpha)}$
for all simple closed curves $\alpha$.  Since any positive power of a Dehn twist about a curve determines that curve, we have that 
$\phi^2(\alpha)=\alpha$ for each $\alpha$.  It follows that either $\phi^2={\rm
Id}$ or $\genus(S)=2$ and $\phi^2$ is the hyperelliptic involution; our assumption that $\genus(S)>2$ rules out the second possibility, so that $\phi={\rm Id}$, and so commutes with 
any element $\psi_2\in\Mod(S)$.  We then pick $\psi_2\in H$ to have infinite 
order.   So we can assume $\phi^2={\rm Id}$ and $\phi\neq {\rm Id}$.

Now any element $\phi\in\Mod^\pm(S)$ of order $2$ is
represented by a homeomorphism $\phi$ of order $2$ (by a
theorem of Fenchel).  We now assume that $g=\genus(S)>2$.  
First suppose that $\Fix(\phi)$ is discrete.  Then $S$ two-fold 
branched covers $S/\langle\phi\rangle$.  Since $g=\genus(S)\geq 3$, the
Riemann-Hurwitz formula easily implies that either
$\genus(S/\langle\phi\rangle)>0$ or that there are at least $4$ branch points
on $S/\langle\phi\rangle$.  Either way, the quotient $S/\langle\phi\rangle$ admits a
self-homeomorphism $\psi$ whose mapping class has infinite order.
After perhaps replacing $\psi$ by a finite power of $\psi$, we
know that $\psi$ lifts to a self-homeomorphism $\psi_2$ of $S$
with the property that, in $\Mod(S)$ we have
$\psi_2\phi=\phi\psi_2$.  By replacing $\psi_2$ with an appropriate power if necessary, we may assume that $\psi_2$ lies in the finite index subgroup $H$.  

If $\Fix(\phi)$ is not discrete then $\phi$ is orientation-reversing and $\Fix(\phi)$ is a union 
of $c>0$ simple closed curves.  If the quotient $S':=S/\langle\phi\rangle$ has $\genus(S')>0$, then $S'$ admits an infinite order self-homeomorphism, which we can then lift as above to obtain $\psi_2$.  If $\genus(S')=0$ 
then $S'$ is planar.  Picking the outermost curve gives $S'$ the structure of a disk with 
$(c-1)$ open disks removed from its interior.  Thus the Euler characteristic 
$\chi(S')=1-(c-1)=2-c$.  Since $S$ is obtained from $S'$ by gluing $2$ copies of $S'$ along its ($\chi=0$) 
boundary, we have $2-2g=\chi(S)=4-2c$ so that $c=g+1$.  Since we are assuming $g>2$, it is clear that $S'$ has an infinite order self-homeomorphism, and we are done as above.
\endproof

\section{Proof of Theorem \ref{theorem:exact}}

By Theorem \ref{theorem:main}, $[\Isom(\Teich(S)):\Mod(S)]<\infty$.  We pass to the index $2$ subgroup 
$\Isom^+(\Teich(S))$ of orientation-preserving isometries of $\Teich(S)$.  Note that 
that any element $f\in \Isom^+(\Teich(S))$ must have $\Fix(f)$ of codimension at least $2$.  Let 
$f\in \Isom^+(\Teich(S))$ with $f\not\in\Mod(S)$ be given.  Ivanov's theorem on commensurations of $\Mod(S)$ mentioned above implies that the conjugation action of $f$ on some characteristic finite index  subgroup $H'$ of $\Mod(S)$ agrees on some finite index subgroup $H\leq H'$ 
with conjugation by some element $\phi\in\Mod(S)$.  By 
composing with $\phi$, we may assume the conjugation action of $f$ on $H$ is trivial, i.e. that 
$f$ centralizes $H$.  As $[\Isom^+(\Teich(S)):H]<\infty$, it must be that $f^n\in H$ for some $n>1$.   
Now consider the exact sequence
\begin{equation}
\label{eq:3}
1\to H\to \langle H,f \rangle \to \langle f\rangle/\langle f^n\rangle\to 1
\end{equation}
As $H$ is centerless (e.g. since it is finite index in $\Mod(S)$ and so contains a pair of independent pseudo-Anosovs) and since the action of $f$ on $H$ is trivial, it follows that (\ref{eq:3}) splits, 
so that $f^n={\rm Id}$. By passing to a power of $f$ if
necessary, we may assume that $f$ has order some prime $p\geq 2$.  Hence $\Fix(f)$ is $\Z/p\Z$-acyclic by Smith theory.  Now $\Fix(f)$ has codimension at least $2$, and so has dimension
at most $2$.  It is also a manifold.  Since $\dim(\Fix(f))\leq
2$, it follows that $\Fix(f)$ is contractible.   But, just as
noted in Case (B) of Step 3 in \S\ref{section:royden} above, $H$ in these cases is not the fundamental
group of a (closed or open) surface; it is also not the
fundamental group of a $1$-manifold by the same argument.  We
thus have a contradiction, so that $f$ must be trivial.
\endproof

\noindent
Dept. of Mathematics, University of Chicago\\
5734 University Ave.\\
Chicago, Il 60637\\
E-mail: farb@math.uchicago.edu, shmuel@math.uchicago.edu
\medskip

\end{document}